\RequirePackage{ifpdf}
\ifpdf 
\documentclass[pdftex]{sigma}
\else
\documentclass{sigma}
\fi

\begin{document}

\allowdisplaybreaks

\renewcommand{\thefootnote}{$\star$}

\renewcommand{\PaperNumber}{066}

\FirstPageHeading

\ShortArticleName{Hochschild Cohomology Theories in White Noise
Analysis}

\ArticleName{Hochschild Cohomology Theories\\ in White Noise
Analysis\footnote{This paper is a contribution to the Special
Issue on Deformation Quantization. The full collection is
available at
\href{http://www.emis.de/journals/SIGMA/Deformation_Quantization.html}{http://www.emis.de/journals/SIGMA/Deformation\_{}Quantization.html}}}

\Author{R\'emi L\'EANDRE}

\AuthorNameForHeading{R. L\'eandre}

\Address{Institut de Math\'ematiques de Bourgogne, Universit\'e de
Bourgogne, 21000, Dijon, France}

\Email{\href{mailto:Remi.leandre@u-bourgogne.fr}{Remi.leandre@u-bourgogne.fr}}

\ArticleDates{Received June 18, 2008, in f\/inal form September
08, 2008; Published online September 27, 2008}

\Abstract{We show that the continuous Hochschild cohomology and
the dif\/ferential Hoch\-schild cohomology of the Hida test
algebra endowed with the normalized Wick product are the same.}

\Keywords{white noise analysis; Hochschild cohomology}

\Classification{53D55; 60H40}

\section{Introduction}

Hochschild cohomology is a basic tool in the deformation theory of
algebras. Gerstenhaber has remarked that in his seminal work  (we
refer to \cite{Gerstenhaber} and references therein for that).
Deformation quantization~\cite{Bayen1,Bayen2} in quantum f\/ield theory leads to some
important problems \cite{Dito1, Dito2, Duetsch,Witten}. Motivated by
that, Dito~\cite{Dito3} has def\/ined the Moyal product on a
Hilbert space. It is easier to work with models of stochastic
analysis although they are similar to models of quantum f\/ield.
In order to illustrate the dif\/ference between these two
theories, we refer to:
\begin{itemize}\itemsep=0pt
\item The paper on Dirichlet forms in inf\/inite dimensions of
Albeverio--Hoegh--Krohn \cite{Albeverio} which used measures on
the space of distributions, the traditional
 space of quantum f\/ield theory.

\item The seminal paper of Malliavin \cite{Malliavin} which used
the traditional Brownian motion and the space of continuous
functions as a topological space.
 This allows  Malliavin to
introduce stochastic dif\/ferential equations in
inf\/inite-dimensional analysis, and to interpret some traditional
tools of quantum f\/ield theory  in stochastic analysis.
\end{itemize}

This remark lead Dito and L\'eandre \cite{DitoLeandre} to
construct of the Moyal product on the Malliavin test algebra on
the Wiener space.

It is very classical in theoretical physics \cite{Di Fransesco}
that the vacuum expectation of some operator algebras  on some
Hilbert space is formally represented by formal path integrals on
the f\/ields. In the case of inf\/inite-dimensional Gaussian
measure, this isomorphism is mathematically well established and
is called the Wiener--It\^o--Segal isomorphism between the Bosonic
Fock space and the $L^2$ of a Gaussian measure. The operator
algebra is the algebra of annihilation and creation operators with
the classical commutation relations. In the case of the classical
Brownian motion~$B_t$ on~$\mathbb{R}$, $B_t$~is identif\/ied with
$A_t+A_t^*$ where $A_t$ is the annihilation operator associated to
$1_{[0,t]}$ and $A_t^*$ the associated creation operator.

Let us give some details on  this identif\/ication~\cite{Meyer4}.
Let $H$ be the Hilbert space of $L^2$ maps from~$\mathbb{R}^+$
into~$\mathbb{R}$. We consider the symmetric tensor product
$\hat{H}^{\otimes n}$ of $H$. It can be realized as the set of
maps $h^n$ from $(\mathbb{R}^{+})^n$ into $\mathbb{R}$ such that
$\int_{(\mathbb{R}^+)^n}\vert h^n(s_1,\dots,s_n)\vert^2ds_1\cdots
ds_n = \Vert h^n \Vert^2 < \infty$. Moreover these maps
$h^n(s_1,\dots,s_n)$ are symmetric in $(s_1,\dots,s_n)$. The
symmetric Fock space~$W_0$ coincides with the set of formal series
$\sum h^n$ such that $\sum n!\Vert h^n\Vert^2 < \infty$. The
annihilation opera\-tors $A_t$ and the creation operators are
densely def\/ined on~$W_0$, mutually adjoint and therefore
closable. To $h^n$ we associate the multiple Wiener chaos
$I^n(h^n)$
\begin{gather*}
I^n(h^n) =\int_{(\mathbb{R}^+)^n}h^n(s_1,\dots,s_n)\delta
B_{s_1}\cdots \delta B_{s_n},
\end{gather*}
where $s \rightarrow B_s$ is the classical Brownian motion on
$\mathbb{R}$. $E_{dP}[\vert I^n(h^n)^2\vert]$for the law of the
Brownian motion $dP$ is $n!\Vert h^n \Vert^2$. $I^n(h^n)$ and
$I^m(h^m)$ are mutually orthogonal on~$L^2(dP)$. If $F$ belongs
to~$L^2(dP)$, $F$ can be written in a unique way $F= \sum
I^n(h^n)$ where $\sum h^n$ belongs to the symmetric Fock
space~$W_0$. This identif\/ication, called chaotic decomposition
of Wiener functionals, realizes an isometry between $L^2(dP)$ and
the symmetric Fock space. $B_t$ can be assimilated to the densely
def\/ined closable operator on $L^2(dP)$
\begin{gather*}
F \rightarrow B_tF.
\end{gather*}
This operator is nothing else but the operator $A_t + A_t^*$ on
the symmetric Fock space.

White noise analysis~\cite{Hida1}  is concerned with the time derivative
of $B_t$ (the white noise) as a distribution  (an element of
$W_{-\infty}$) acting on some weighted Fock space $W_{-\infty}$
(we refer to \cite{Berezanskii,Hida2,Obata2} for textbooks on
white noise analysis). Let us recall namely that the Brownian
motion is only continuous! The theory of Hida distribution leads
to new insight in stochastic analysis.

One of the main points of interest in the white noise analysis is
that we can compute the elements of $L(W_{\infty-},W_{-\infty})$
\cite{Ji1,Ji2, Obata2} in terms of kernels.
 We refer to \cite{Ji2, Maassen, Meyer1}
for a well established theory of kernels on the Fock space.
 This theory was motivated by the heuristic constructions of quantum f\/ield theory \cite{Berezin1, Glimm, Haag}. Elements
of $L(W_{\infty-}, W_{-\infty})$ can be computed in a sum of
multiple integrals of the elementary creation and annihilation
operators.

This theorem  plays the same role as the theorem of
Pinczon~\cite{Pinczon1,Pinczon2}: the  operators acting on
$\mathbb{C}(x_1,\dots,x_d)$, the complex polynomial algebra on
$\mathbb{R}^d$, are series of dif\/ferential operators with
polynomial components. This theorem of Pinczon allowed  Nadaud
\cite{Nadaud1, Nadaud2} to show that the continuous Hochschild
cohomology on $C^\infty(\mathbb{R}^d)$ is equal to the
dif\/ferential Hochschild cohomology of the same algebra
 (we refer to  papers of Connes \cite{Connes} and Pf\/laum \cite{Pflaum} for other proofs).

In the framework of white noise analysis we have an analogous
theorem to the theorem of Pinczon~\cite{Pinczon1}. Therefore we
can
 repeat in this framework the proof of Nadaud. We show that the
continuous Hochschild cohomology~\cite{Leandre1} of the Hida Fock space (we
consider {\it series} of kernels) is equal to the dif\/ferential
Hochschild cohomology (we consider {\it finite} sums of kernels).

In the f\/irst part of this work, we recall the theorem of Obata
which computes the operators on the Hida Fock space: Obata
considers standard creation operators and standard annihilation
operators. We extend this theorem in the second part to continuous
multilinear operators on the Hida algebra, endowed with the
normalized Wick product. This Hida test algebra was used by
L\'eandre \cite{Leandre2,Leandre3} in order to def\/ine some star
products in white noise analysis.

We refer to the review paper of L\'eandre for deformation
quantization in inf\/inite-dimensional analysis~\cite{Leandre4}.

\section{A brief review on Obata's theorem}

We consider the Hilbert space $H= L^2(\mathbb{R},dt)$. We consider
the operator $\Delta = 1+t^2-d^2/dt^2$. It  has eigenvalues $\mu_j
= (2j+2)$ associated to the normalized eigenvectors $e_j$, $j \geq
0$. We consider the Hilbert space $H_k$ of series $f = \sum
\lambda_j e_j$ such that
\begin{gather*}
\Vert f \Vert_k^2 = \sum \vert \lambda_j \vert^2 \mu_j^{2k}<
\infty.
\end{gather*}
It is the Sobolev space associated to $\Delta^k$. Since $\mu_j
>1$,  $H_k \subseteq H_{k'}$ if $k>k'$, and the system of Sobolev
norms $\Vert \cdot \Vert_k$ increases when $k$ increases.
Therefore we can def\/ine the test functional space $H_{\infty-}$
of functions $f$ such that all norms $\Vert f\Vert_k < \infty$.
 A functional $F$ with values in this space is continuous if and
only if it is continuous for all Sobolev norms $\Vert  \cdot
\Vert_k$, $k>0$,
\begin{gather*}
H_{\infty-}=\cap_{k\geq 0}H_k.
\end{gather*}
The topological dual of $H_{\infty-}$ is the space of Schwartz
distributions:
\begin{gather*}
H_{-\infty} = \cup_{k<0}H_k.
 \end{gather*}
$\sigma$ is called a distribution if the following condition holds:  let $f$ be in $H_{\infty}$. For some $k>0$, there exists
$C_k$ such that for all $f\in H_{\infty-}$, $\vert \langle
\sigma,f\rangle\vert \leq C_k\Vert f\Vert_k$ Therefore we get a
Gel'fand triple
\begin{gather*}
H_{\infty-}\subseteq H \subseteq H_{-\infty}.
\end{gather*}
We complexify all these spaces (We take the same notation). It is
important to complexify these spaces to apply Potthof\/f--Streit
theorem~\cite{Streit}.

Let $A = ((i_1, r_1),\dots,(i_n,r_n))$ where $i_1 <
i_2<\cdots<i_n$ and $r_i > 0$. We put
\begin{gather}\label{eq7}
 \vert A \vert = \sum r_i
 \end{gather}
and
\begin{gather*}e_A = \hat{\otimes}e_{i_j}^{r_j},\end{gather*}
where we consider a normalized symmetric tensor product. We
introduce the Hida weight
\begin{gather*}
\Vert A\Vert = \prod_{(i_j,r_j) \in A} (2i_j+2)^{r_j}.
\end{gather*}
We consider the weighted Fock space $W_k$ of series
\begin{gather*}
\phi = \sum_A \lambda_A e_A
\end{gather*}
such that
\begin{gather*}
\Vert \phi \Vert_k = \sum_A \vert \lambda_A\vert^2 \Vert
A\Vert^{2k} \vert A \vert! < \infty
\end{gather*}
($\lambda_A$ is complex). These systems of norms increase when $k$
increases.

We consider
\begin{gather*}W_{\infty-} = \cap_{k>0}W_k\end{gather*}
endowed with the projective topology and its topological
dual (called the space of Hida distributions)
\begin{gather*}W_{-\infty} = \cup_{k<0}W_k\end{gather*}
endowed with the inductive topology. We get a Gel'fand triple
\begin{gather*}W_{\infty-} \subseteq W_0 \subseteq W_{-\infty},
\end{gather*}
$W_0$ is the classical Fock space of quantum physics.

We consider $\xi \in H_{\infty-}$ and the classical
coherent vector
\begin{gather*}
\phi_\xi = \sum_{n\in \mathbb{N}} {\xi^{\otimes n}\over n!},
\end{gather*}
$\phi_\xi$ belongs to the Hida test functional space
$W_{\infty-}$.
\begin{definition}
Let $\Xi$ belong to $L(W_{\infty-}, W_{-\infty})$. Its symbol is
the function $\hat{\Xi}$ from $H_{\infty-}\times H_{\infty-}$ into
$\mathbb{C}$ def\/ined by
\begin{gather*}
\hat{\Xi}(\xi, \eta) = \langle \Xi
\phi_{\xi},\phi_{\eta}\rangle_0.
\end{gather*}
\end{definition}

If $\Xi$ belongs to $L(W_{\infty-}, W_{-\infty})$, its symbol
satisf\/ies clearly the following properties:
\begin{enumerate}\itemsep=0pt
\item[(P${}_1)$] For any $\xi_1, \xi_2, \eta_1, \eta_2 \in
H_{\infty-}$, the function
\begin{gather*} (z,w) \rightarrow \hat{\Xi}(z\xi_1+\xi_2, w\eta_1 + \eta_2)\end{gather*}
is an entire holomorphic function on $\mathbb{C} \times \mathbb{C}$.

\item[(P${}_2$)] There exists a constant $K$ and a constant $k>0$
such that
\begin{gather*}
\vert \hat{\Xi}(\xi, \eta)\vert^2 \leq K \exp\big[\Vert
\xi\Vert_k^2\big]\exp\big[\Vert \eta\Vert_k^2\big].
\end{gather*}
\end{enumerate}

The converse of this theorem also holds. It is a result
of Obata~\cite{Obata1} which generalizes the theorem of
Potthof\/f--Streit
 characterizing distribution in
white noise analysis~\cite{Streit}. If a~function $\hat{\Xi}$ from
$H_{\infty-}\times H_{\infty-}$ into $\mathbb{C}$ satisf\/ies
(P${}_1$) and (P${}_2$), it is the symbol
 of an element of $L(W_{\infty-}, W_{-\infty})$. Its continuity norms can be estimated in a universal way linearly in the data of (P${}_2$).

This characterization theorem allows Obata to show that an element
$\Xi$ of $L(W_{\infty-}, W_{-\infty})$ can be decomposed into a
sum
\begin{gather*} \Xi = \sum_{l,m}\Xi_{l,m}(k_{l,m}),
\end{gather*}
where $\Xi_{l,m}(k_{l,m})$ is def\/ined by the following
considerations.

Let $A$ be as given above \eqref{eq7}. Let $a_{i}^*$ be the
standard creation operator
\begin{gather*}
a_i^*e_A = c(r_i)e_{A^i},
\end{gather*}
where
\begin{gather*}
A^i = ((i_1,r_1),\dots,(i_l,r_l), (i,r_i+1),
(i_{l+1},r_{l+1}),\dots,(i_n,r_n))
\end{gather*}
if $i_l< i<i_{l+1}$. If $i$ does not appear in $A$, we put $r_i =
1$. We consider the standard annihilation operator $e_i$ def\/ined
\begin{gather*}a_ie_A = 0\end{gather*}
if $i$ does not appear in $A$ and equals to $c'(r_i)e_{A_i}$ where
\begin{gather*}A_i = ((i_1,r_i),\dots, (i_l,r_l), (i,r_i-1),\dots,(i_n,r_n)).
\end{gather*}
The constants $c(r_i)$ and $c'(r_i)$ are computed
in~\cite{Meyer4}. Their choice is motivated by the use of Hermite
polynomial on the associated Gaussian space by the
Wiener--It\^o--Segal isomorphism: the role of Hermite polynomial
in inf\/inite dimensions is played by the theory of chaos
decomposition through the theory of multiple Wiener integrals. The
annihilation operator $a_i$ corresponds to the stochastic
derivative in the direction of $e_i$ on the corresponding Wiener
space. $a_i^*$ is its adjoint obtained by integrating by parts on
the Wiener space. We consider the operator $\Xi_{I,J}$
\begin{gather*} \phi \rightarrow a_{i_1}^*\cdots a_{i_l}^*a_{j_1}\cdots a_{j_m}\phi = a^*_Ia_J\phi.\end{gather*}
It belongs to $L(W_{\infty-},W_{-\infty})$ and its symbol
is~\cite{Obata1}
\begin{gather*}
\exp[\langle \xi,\eta\rangle _0]\prod_{j_k\in J} \langle
e_{j_k},\xi\rangle_0\prod_{i_k\in I}\langle e_{i_k},\eta\rangle_0.
\end{gather*}
Therefore we can consider
\begin{gather*} \Xi_{l,m}(k_{l,m}) = \sum_{\vert I \vert=l, \vert J\vert = m}\lambda_{I,J}a^*_Ia_J,
\end{gather*}
where
\begin{gather*}\sum \vert \lambda_{I,J}\vert^2\Vert I \Vert^{-k}\Vert J\Vert^{-k} < \infty\end{gather*}
for some $k>0$. $\Xi_{l,m} = \sum\lambda_{I,J}e_I\otimes e_J$
def\/ines an element of $H_{-\infty}^{\otimes (l+m)}$. $\Xi_{I,J}$
can be extended by linearty to
\begin{gather*}
\sum\lambda_{I,J}\Xi_{I,J} = \Xi_{l,m}(k_{l,m}).
\end{gather*}
$\Xi_{l,m}(k_{l,m})$ belongs to $L(W_{\infty-}, W_{-\infty})$ if $k_{l,m}$ belongs to $H_{-\infty}^{\otimes (l+m)}$.
 This last space is $\cup_{k>0} H_{-k}^{\otimes (l+m)}$
 endowed with the inductive topology. $\Xi_{l,m}(k_{l,m})$ belongs to
$L(W_{\infty-}, W_{\infty-})$ if $k_{l,m}$ belongs to
$H_{\infty-}^{\otimes l}\otimes H_{-\infty}^{\otimes m}$.  This
means that there exists $k$
 such that for all $k'$
\begin{gather*}\sum \vert \lambda_{I,J}\vert^2\Vert I \Vert^{k'}\Vert J\Vert^{-k} < \infty\end{gather*}
for some $k>0$. The symbol of $\Xi_{l,m}(k_{l,m})$ satisf\/ies
\begin{gather*}
\hat{\Xi}_{l,m}(k_{l,m})(\eta,\xi) = \langle k_{l,m},
\eta^{\otimes l}\otimes \xi^{\otimes m}\rangle \exp[\langle \xi,
\eta\rangle _0].
\end{gather*}

Following the heuristic notation of quantum f\/ield  theory
\cite{Berezin1,Berezin2,Glimm,Haag,Meyer3}, the operator $\Xi_{l,m}(k_{l,m})$ can
be written as
\begin{gather*}
\Xi_{l,m}(k_{l,m}) = \int_{\mathbb{R}^{l+m}}k(s_1,\dots,s_l,
t_1,\dots,t_m)a^*_{s_1}\cdots a^*_{s_l}a_{t_1}\cdots
a_{t_m}ds_1\cdots ds_ldt_1\cdots dt_m.\end{gather*}
 The ``elementary'' creation operators $a^*_s$ and the ``elementary'' annihilation operators $a_t$ satisfy  the canonical commutation relations
\begin{gather*} [a^*_s,a^*_t] = [a_s, a_t] = 0,\qquad [a^*_s,a_t] = \delta(s-t),
\end{gather*}
where $\delta(\cdot)$ is the Dirac function in~0.

\section{Fock expansion of continuous multilinear operators}

We are motivated in this work by the Hochschild cohomology in
white noise analysis. For that, we require that $W_{\infty-}$ is
an algebra.
 In order to be self-consistent we will take the model of \cite{Leandre3} or \cite{Leandrerogers}.

We will take the normalized Wick product
\begin{gather*} :e_A.e_B:  = e_{A\cup B}, \end{gather*}
where $A \cup B$ is obtained by concatenating the indices and
adding the length of these when the same appears twice.

$W_{\infty-}$ is not the same space as before. We consider another
Hida Fock space. $W_{k,C}$ is the space of $\phi = \sum_A
\lambda_A e_A$ such that
\begin{gather*}\Vert \phi\Vert_{k,C}^2 = \sum_A \vert \lambda_A\vert^2C^{2\vert A\vert}\Vert A\Vert^{2k} \vert A\vert! < \infty. \end{gather*}
$W_{k,C}$ can be identif\/ied with the Bosonic Fock space
associated to the Hilbert Sobolev space associated to the operator
$C\Delta^k$. $W_{\infty-}$ is the intersection of $W_{k,C}$,
$k>0$, $C>0$. This space is endowed with the projective topology.

By a small improvement of \cite{Leandre3} and
\cite{Leandrerogers}, we get:
\begin{theorem}
$W_{\infty-}$ is a topological algebra for the normalized Wick
product.
\end{theorem}

\begin{proof} The only new ingredient in the proof of \cite{Leandre3,Leandrerogers} is that
\begin{gather*}
\vert A \cup B\vert! \leq 2^{\vert A \vert + \vert B\vert}\vert
A\vert! \vert B\vert!.
\end{gather*}
Let us give some details. Let us consider
\begin{gather*} \phi^1 = \sum_A \lambda^1_A e_A,\qquad
 \phi^2 = \sum_A \lambda^2_A e_A.\end{gather*}
We have
\begin{gather*} :\phi^1\cdot \phi^2: = \sum_A \mu_A e_A,\qquad \mbox{where}\qquad
 \mu_A = \sum_{B\cup D = A}\lambda^1_B \lambda^2_D. \end{gather*}
There are at most $2^{\vert A\vert}$ terms in the previous sum. By
Jensen inequality
\begin{gather*} \vert \mu_A\vert^2 \leq C_1^{\vert A\vert} \sum_{B\cup D = A}\vert \lambda_B^1\vert^2 \vert \lambda_D^2\vert^2.\end{gather*}
Therefore
\begin{gather*}\Vert:\phi^1\cdot \phi^2:\Vert^2_{k,C} \leq \sum_A(C_1C)^{2\vert A\vert}C^{2\vert A\vert}\Vert A\Vert^{2k}\vert A\vert!
\sum_{B\cup D = A}\vert \lambda_B^1\vert^2 \vert
\lambda_D^2\vert^2.\end{gather*} But
\begin{gather*} \Vert A\Vert^{2k} \leq \Vert B \Vert^{2k} \Vert D\Vert^{2k}\qquad \mbox{and}\qquad
 \vert A\vert! \leq 2^{\vert B\vert + \vert D\vert}\vert B\vert! \vert D\vert!\end{gather*}
if the concatenation of $B$ and $D$ equals $A$. Therefore, for
some $C_1$
\begin{gather*}\Vert :\phi^1\cdot \phi^2:\Vert^2_{k,C} \leq \sum_A \sum_{B \cup D = A}(C_1C)^{\vert B\vert + \vert D\vert}
\Vert B \Vert^{2k} \Vert D\Vert^{2k} \vert B\vert! \vert D\vert!
\leq   \Vert \phi^1\Vert^2_{k, C_1C} \Vert \phi^2\Vert^2_{k,
C_1C}.\!\!\!
\end{gather*}
This shows the result.
\end{proof}

Let $L(W^n_{\infty-},W_{\infty-})$ be the space of $n$-multilinear
continuous applications from $W_{\infty-}$ in\-to~$W_{\infty-}$.

\begin{definition}
The symbol $\hat{\Xi}$ of an element $\Xi$ of
$L(W^n_{\infty-},W_{\infty-})$
 is the map from $H_{\infty-}^n \times H_{\infty-}$ into $\mathbb{C}$ def\/ined by
\begin{gather*}
\xi^1\times \xi^2\times\cdots \times \xi^n \times \eta \rightarrow
\langle \Xi(\phi_{\xi^1},\dots,\phi_{\xi^n}), \phi_\eta\rangle_0 =
\hat{\Xi}(\xi^1,\dots,\xi^n,\eta).
\end{gather*}
\end{definition}

$\Xi$ belongs to $L(W_{\infty-}^n, W_{\infty-})$ if for any
$(k,C)$, there exists $(k_1, C_1, K_1)$ such  that
\begin{gather*}
\Vert \Xi(\phi^1,\dots,\phi^n)\Vert_{k,C}\leq K_2 \prod
\Vert\phi^i\Vert_{k_2,C_2}.
\end{gather*}

If $\Xi$ belongs to $L(W^n_{\infty-}, W_{\infty-})$, its symbol
satisf\/ies clearly  the following properties:
\begin{enumerate}\itemsep=0pt
\item[(O${}_1$)] For any elements $\xi^1_1,\dots,\xi^n_1$,
$\xi^1_2,\dots,\xi^n_2$, $\eta_1$, $\eta_2$ of $H_{\infty-}$, the
map
\begin{gather*}(z_1,\dots,z_n, w) \rightarrow \hat{\Xi}(z_1\xi^1_1 + \xi_2^1,\dots, z_n\xi_1^n+\xi_2^n, w\eta_1 + \eta_2)\end{gather*}
is an entire holomorphic map from $\mathbb{C}^n \times
\mathbb{C}$ into $\mathbb{C}$.

\item[(O${}_2$)] For all $k>0$, $K>0$, there exists numbers $C$,
$k_1>k$, $K_1>0$ such that
\begin{gather*}\vert \hat{\Xi}(\xi^1,\dots,\xi^n, \eta)\vert^2 \leq C\exp\left[K_1\sum_{i=1}^n\Vert \xi^i\Vert^2_{k_1}+ K\Vert \eta\Vert_{-k}^2\right]. \end{gather*}
\end{enumerate}

We prove
the converse of this result. It is a small improvement of the
theorem of Ji and  Obata~\cite{Ji1}.

\begin{theorem}
If a function $\hat{\Xi}$ from $\mathbb{C}^n \times \mathbb{C}$
into $\mathbb{C}$ satisfies to {\rm (O${}_1$)} and {\rm
(O${}_2$)},
 it is the symbol
 of an element $\Xi$ of $L(W_{\infty-}^n, W_{\infty-})$. The different modulus of continuity can be estimated in terms of the data
in {\rm (O${}_2$)}.
\end{theorem}

\begin{proof}
It is an adaptation of the proof of  a result of
Obata~\cite{Obata1}, the result which was generalizing
Potthof\/f--Streit theorem. We omit all the details. This
classical Potthof\/f--Streit theorem is the following. Let $\Phi$
in $W_{-\infty}$ be the topological dual of $W_{\infty-}$. We
def\/ine its $S$-transform as follows
\begin{gather*}S(\xi) = \langle \Phi, \phi_\xi\rangle_0\end{gather*}
for $\xi \in H_{-\infty}$. The $S$-transform of $\Phi$ satisf\/ies
the following properties:
\begin{enumerate}\itemsep=0pt
\item[i)] the function
$z \rightarrow S(z\xi_1 + \xi_2)$
is entire holomorphic;

\item[ii)] for some $K_1>0$, $K_2 >0$ and some $k\in \mathbb{R}$
\[
\vert S(\xi)\vert^2 \leq K_1 \exp[K_2\Vert \xi\Vert^2_k].
\]
\end{enumerate}

Potthof\/f--Streit theorem states the opposite \cite[Lemma
3.2]{Ji1}: if a function $S$ from $H_{\infty-}$ into~$\mathbb{C}$
satisf\/ies i) and ii),
 it is the $S$-transform
of a distribution~$\Phi$. Moreover,  there exists $C$ depending
only of $K_2$ such that
\begin{gather*}
\Vert \Phi\Vert^2_{-k-r, C} \leq K_1\end{gather*} for all $r>0$.

From this theorem, we deduce that there exists a distribution
$\Phi_{\xi^1,\dots,\xi^{n-1}, \eta}$ such that
\begin{gather*}
\hat{\Xi}(\xi^1,\dots,\xi^n,\eta) = \langle
\Phi_{\xi^1,\dots,\xi^{n-1},\eta},\phi_{\xi}\rangle
_0.\end{gather*} Moreover there exists $C$ independent of $\eta$,
$\xi^1,\dots,\xi^{n-1}$ such that
\begin{gather*}
\Vert \Phi_{\xi^1,\dots,\xi^{n-1},\eta}\Vert^2_{-k_1-r,C} \leq
K_2\exp\left[K_1\sum_{i=1}^{n-1}\Vert \xi^i\Vert^2_{k_1}+K\Vert
\eta\Vert^2_{-k}\right].
\end{gather*}
If $\phi$ belongs to $W_{\infty-}$ we put
\begin{gather*}
G_\phi(\xi^1,\dots,\xi^{n-1}, \eta) = \langle
\Phi_{\xi^1,\dots,\xi^{n-1}, \eta},\phi\rangle_0.
\end{gather*}
We have for some $K_2$, $C_2$, $K_1$, $k_1$, $k_2$ depending only on the previous datas that
\begin{gather*}
\Vert G_\phi(\xi^1,\dots,\xi^{n-1}, \eta)\Vert^2 \leq
K_2\Vert\phi\Vert_{k_2,C_2}^2 \exp\left[K_1\sum_{i=1}^{n-1}\Vert
\xi^i\Vert_{k_1}^2+K\Vert \eta\Vert_{-k}^2\right].
\end{gather*}
The two properties (O${}_1$) and (O${}_2$) are satisf\/ied at the
step $n-1$. By induction, we deduce that
\begin{gather*}
G_\phi(\xi^1,\dots,\xi^{n-1},\eta) =
\hat{\Xi}_\phi(\xi^1,\dots,\xi^{n-1},\eta).
\end{gather*}
Moreover, we get that
\begin{gather*}G_\phi(\xi^1,\dots,\xi^{n-1},\eta) = \langle \Xi_\phi(\phi_{\xi^1},\dots,\phi_{\xi^{n-1}}), \phi_\eta\rangle_0,
\end{gather*}
where $\Xi_\phi$ is an element of $L(W_{\infty-}^{n-1},
W_{\infty-})$ depending linearly and continuously from $\phi \in
W_{\infty-}$. We put
\begin{gather*} \Xi(\phi^1,\dots,\phi^{n-1},\phi) = \Xi_\phi(\phi^1,\dots,\phi^{n-1}).\end{gather*}
It remains to prove the result for $n = 1$. It is a small
improvement of the proof of the result of Ji and
Obata~\cite{Ji1}. Let us give some details.

By using Potthof\/f--Streit theorem, we deduce that there is a
distribution $\Phi_\eta$ such that
\begin{gather*}\hat{\Xi}(\xi, \eta) = \langle \Phi_\eta, \phi_\xi\rangle_0.
\end{gather*}
Moreover there exists $C$ independent of $\eta$ such that
\begin{gather*}\Vert \Phi_\eta\Vert^2_{-k_1-r,C}\leq K_2\exp[K\Vert\eta\Vert^2_{-k}]\end{gather*}
If $\phi$ belongs to $W_{\infty-}$, we set
\begin{gather*}G_\phi(\eta) = \langle \Phi_\eta,\phi\rangle_0.\end{gather*}
We have
\begin{gather*}\vert G_\phi(\eta)\vert^2 \leq K_3\Vert \phi\vert_{k_2,C_2}^2\exp[K\Vert\eta\Vert_{-k}^2]\end{gather*}
for some $k_2 > 0$.  We apply Potthof\/f--Streit theorem
(see~\cite[Lemma 3.3]{Ji1}). There exists an element $\Xi(\phi)$
of $W_{k-r, C}$ where $k>0$ depending continuously of $\phi$ such
that
\begin{gather*}G_\phi(\eta) = \langle \Xi(\phi), \phi_\eta\rangle. \end{gather*}
We have clearly
\begin{gather*}\langle \Xi(\phi_\xi), \phi_\eta\rangle  = \hat{\Xi}(\xi, \eta).\end{gather*}
This shows the result.
\end{proof}

The following statements follow closely \cite[Appendix]{Chung}.

Let $\Xi$ be an element of $L(W^n_{\infty-}, W_{\infty-})$. Let
$\hat{\Xi}$ be its symbol. We put:
\begin{gather*}
\Psi(\xi^1,\dots,\xi^n, \eta) = \exp\left[-\sum_{i=1}^n\langle
\xi_i, \eta\rangle_0\right]\hat{\Xi}(\xi^1,\dots,\xi^n, \eta).
\end{gather*}

Clearly $\Psi$ satisf\/ies to (O${}_1$) and (O${}_2$). We put
\begin{gather*}
\psi(z_1^1,\dots,z_{m_1}^1,z_1^2,\dots,z_{m_2}^2,\dots,z_1^n,\dots,z_{m_n}^n,w_1,\dots,w_l) \nonumber\\
\qquad
{}=\Psi(z_1^1\xi_1^1+\cdots+z_{m_1}^1\xi^1_{m_1},\dots,z_1^n\xi_1^n+\cdots
+z^n_{m_n}\xi^n_{m_n}, w_1\eta_1+\cdots+w_l\eta_l).
\end{gather*}
We put $M = (m_1,\dots,m_n)$ and
\begin{gather*}
K_{l,M}(\xi_1^1,\dots,\xi_{m_1}^1,\dots,\xi_1^n,\dots,\xi_{m_n}^n,\eta_1,\dots,\eta_l) \nonumber \\
\qquad{}={1\over l!m_1!\cdots m_n!}{\partial^{l+\sum m_i}\over
\partial z_1^1\cdots \partial z_{m_1}^1\cdots \partial
 z_1^n \cdots \partial z^n_{m_n}\partial w_1\cdots \partial w_l}\psi(0,0,\dots,0).
 \end{gather*}
$K_{L,M}$ is an $l + \sum m_i$ multilinear map.

Since $\psi$ is holomorphic, we have a Cauchy type representation
of the considered expression
$K_{l,M}(\xi_1^1,\dots,\xi_{m_1}^1,\dots,\xi_1^n,\dots,\xi_{m_n}^n,\eta_1,\dots,\eta_l)$
\begin{gather*}
K_{l,M}(\xi_1^1,\dots,\xi_{m_1}^1,\dots,\xi_1^n,\dots,\xi_{m_n}^n,\eta_1,\dots,\eta_l)
= {1\over l!m_1!\cdots m_n!}\prod_{j=1}^n \prod_{i=1}^{m_j}{1\over
2\pi}\int_{\vert z_i^j\vert =r_i^j}{\vert dz_i^j\vert \over
(z_i^j)^2} \nonumber\\ \qquad{}\times \prod_{k=1}^l\int_{\vert
w_k\vert = s_k}{\vert dw_k\vert\over w_k^2}
\psi(z_1^1,\dots,z_{m_1}^1,\dots,z_1^n,\dots,z_{m_n}^n,w_1,\dots,w_l).
\end{gather*}
 We deduce from (O)${}_2$ a bound of $K_{l,M}$ of the type (D.5) in \cite{Chung}
\begin{gather*} \vert K_{l,M}(\xi_1^1,\dots,\xi_{m_1}^1,\dots,\xi_1^n,\dots,\xi_{m_n}^n,\eta_1,\dots,\eta_l)\vert \leq
{C\over l!\prod m_i!}  {1\over r_1^1\cdots r_{m_1}^1 \cdots r_1^n\cdots r_{m_n}^ns_1\cdots s_l}\nonumber\\
\qquad{}\times \exp\left[K_1\left\{\sum_{i,j}r_i^j\Vert
\xi_i^j\Vert_{k_1}\right\}^2\right]
\exp\left[K\left\{\sum_js_j\Vert \eta_j\Vert_{-k}\right\}^2\right]
\end{gather*}
for some $k_1 > k$ and $K$, $k>0$ and some $K_1 > 0$.

According to \cite[(D.6)]{Chung}, we choose
\begin{gather*}r_i^j = {R\over C m_j\Vert \xi_i^j\Vert_{k_1}},\qquad
s_j = {S\over Cl \Vert \eta_j\Vert_{-k}}\end{gather*} and we
deduce a bound of $K_{l,M}$ in
\begin{gather}\label{eq74}
{C\over l!\prod m_i!}\prod \left({Cm_i\over R}\right)^{m_i}
\left({Cl\over S}\right)^l\prod\Vert \xi_i^j\Vert_{k_1} \prod
\Vert \eta_i\Vert_{-k}\exp\big[KR^2\big]\exp\big[KS^2\big].
\end{gather}
Clearly, $K_{l,M}$ is a multilinear application in $\xi_i^j$,
$w_k$. By \eqref{eq74}, $K_{l,M}$ is continuous. Therefore
$K_{l,M}$ can be identif\/ied with an element of
$H_{\infty-}^{\otimes l}\otimes H_{-\infty}^{\otimes \sum m_i}$.
We consider
 \begin{gather*}
\hat{\Xi}_{l,M}(\xi^1,\dots,\xi^n, \eta) =
K_{l,M}(\xi^1,\dots,\xi^1,\xi^2,\dots,\xi^2,\dots,\xi^n,\dots,\xi^n,\eta,\dots,\eta),
\end{gather*}
where $\xi_i$ is taken $m_i$ times and $\eta$ $l$ times.

By holomorphy,
\begin{gather*}
\hat{\Xi} = \sum \hat{\Xi}_{l,M}\exp\left[\sum_{i=1}^n\langle
\xi^i,\eta\rangle_0\right].
\end{gather*}
and the series converges in the sense of (O${}_1$) and (O${}_2$).
Only the second statement presents some dif\/f\/iculties. We
remark for that by \cite[page~557]{Chung} ($\alpha(n) = 1$,
$G_\alpha(s) = \exp[s]$)
\begin{gather*}
\inf_{s>0}\exp[s]s^{-n} \leq Cn!n^{-2n}.
\end{gather*} We deduce a bound analog to the bound (D.7) in \cite{Chung}
\begin{gather*}
\vert \hat{\Xi}_{l,M}\vert^2 \leq  {1 \over (l^l \prod m_i^{m_i})}
C^l \prod C^{m_i} \exp\left[D\sum \Vert \xi^i\Vert^2_{k_1}+
D_1\Vert \eta\Vert_{-k}^2\right].
\end{gather*}
$x^n\exp[-D_1x^2]$ has a bound in $\exp[-C_1n]C^nn^{n/2}$. If
$D_1$ is large, $C$ can be chosen very small and~$C_1$ very large.
We deduce the following bound
\begin{gather*}
\vert \hat{\Xi}_{l,M}\vert^2 \leq C^l C^{\sum
m_i}\exp\left[-C_1l-C_1\sum m_i\right] \exp\left[D_2\sum \Vert
\xi^i\Vert^2_{k_1}+ D_2\Vert \eta\Vert_{-k}^2\right].
\end{gather*} We remark if $D_2$ is very large that $C$ can be chosen very small and $C_1$ can be chosen very large.
We remark if $C_1$ is large
\begin{gather*}
\sum_{l,M}\exp\left[-C_1l-C_1\sum m_i\right] < \infty
\end{gather*}
in order to see that the series $\Xi_{l,M}$ converges in
$L(W_{\infty-}^n, W_{\infty-})$.

\begin{definition}
The series $\sum_{l,M}\Xi_{l,M} = \Xi$ is called the Fock
expansion of the element $\Xi$ belonging to $L(W_{\infty-}^n,
W_{\infty-})$.
\end{definition}

\section{Isomorphism of Hochschild cohomology theories}

In this part, we prove the main theorem of this work.
\begin{lemma}If $\xi$ belongs to $H_{\infty-}$,
\begin{gather*}\xi^{\otimes n} = :\xi:^n.
\end{gather*}\end{lemma}

\begin{proof} We put $\xi = \sum \lambda_i e_i$ such that
\begin{gather*}
\xi^{\otimes n} = \sum_{i_1,\dots,i_n}\lambda_{i_1}\cdots
\lambda_{i_n}e_{i_1}\otimes_0\cdots \otimes_0e_{i_n},
\end{gather*}
where $\otimes_0$ denotes the traditional tensor product. By
regrouping various element, we deduce that
\begin{gather*}
\xi^{\otimes n} = \sum_{\substack{i<1<\cdots <i_r; \
n_1\not=0,\dots,n_r\not= 0;\\ n_1+\cdots+n_r =
n}}\lambda_{i_1}^{n_1}\cdots \lambda_{i_r}^{n_r}{n!\over
n_1!\cdots n_r!}
e_{i_1}^{\otimes n_1}\hat{\otimes}\cdots \hat{\otimes} e_{i_r}^{n_r}\nonumber \\
\phantom{\xi^{\otimes n}}{} =\sum_{\substack{i<1<\cdots <i_r; \
n_1\not=0,\dots,n_r\not= 0;\\ n_1+\cdots+n_r =
n}}\lambda_{i_1}^{n_1}\cdots \lambda_{i_r}^{n_r}{n!\over
n_1!\cdots n_r!} :e_{i_1}:^{\otimes n_1}\cdots :e_{i_r}:^{n_r} =
:\xi:^n.
\end{gather*}
This shows the result. \end{proof}

\begin{corollary}If $\xi_1 \in H_{\infty-}$ and if $\xi_2 \in H_{\infty-}$
\begin{gather*}\phi_{\xi_1 + \xi_2} = :\phi_{\xi_1}\phi_{\xi_2}:.
\end{gather*}
\end{corollary}

\begin{definition}
Let $\Xi$ belong to $L(W_{\infty-}^r, W_{\infty-})$. Its
Hochschild coboundary $\delta^r$ is def\/ined as follows:
\begin{gather*}
\delta^r\Xi(\phi^1,\dots,\phi^{r+1}) = :\phi^1\Xi(\phi^2,\dots,\phi^r) \nonumber\\
\qquad{}
+\sum_{i=1}^r(-1)^i\Xi(\phi^1,\dots,:\phi^i\phi^{i+1}:,\dots,\phi^{r+1})
 + (-1)^{r+1}:\Xi(\phi^1,\dots,\phi^r)\phi^{r+1}:.
 \end{gather*}
\end{definition}

Classically $\delta^{r+1}\delta^r = 0$.
\begin{definition}We say that an element $\Xi$ of $L(W_{\infty-}^r, W_{\infty-})$ is a homogeneous polydif\/ferential
operator of order $(l,m)$ if its symbol
$\hat{\Xi}(\xi^1,\dots,\xi^r,\eta)$ is equal to
\begin{gather}\label{eq86}
\Psi(\xi^1,\dots,\xi^r,\eta)\exp\left[\sum\langle
\xi^i,\eta\rangle\right],
\end{gather}
where $\Psi$ is a homogeneous polynomial in the  $\xi^i$ of degree
$m$ and in $\eta$ of degree~$l$.
\end{definition}

\begin{proposition}If $\Xi$ is an $r$-polydifferential operator of degree $(l,m)$, $\delta^r\Xi$ is an $(r+1)$-poly\-dif\/fe\-ren\-tial operator of degree $(l,m)$.
\end{proposition}

\begin{proof}Since $:\phi_{\xi^1}\phi_{\xi^2}: = \phi_{\xi^1+\xi^2}$, the only problem is to show that
\begin{gather*}
(\phi^1,\dots,\phi^{r+1}) \rightarrow :{\phi^1}
\Xi(\phi^2,\dots,\phi^{r+1}):
\end{gather*}
is still a polydif\/ferential operator of degree $(l,m)$.

Let $\eta \in H_{\infty-}$ be such that $\Vert \eta\Vert_0 = 1$.
Let us compute
\begin{gather*}
\langle \Xi(\phi_{\xi^2},\dots,\phi_{\xi^{r+1}}),
\phi_{\lambda\eta})\rangle_0 =\Psi(\xi^2,\dots,\xi^{r+1},\lambda
\eta) \exp\left[\lambda\sum_{i=2}^{r+1}\langle
\xi^i,\eta\rangle\right].
\end{gather*}
If we compute the component of
$\Xi(\phi_{\xi^2},\dots,\phi_{\xi^{r+1}}$) along ${\eta^{\otimes
n}\over n!}$, it is the element of degree $n$ in the expansion in
$\lambda$ of the~\eqref{eq86}. Since $\Psi$ is homogeneous of
degree $l$ in $\eta$, the term of degree $n$ in the expansion in
$\lambda$ of~\eqref{eq86} is
\begin{gather*}{\sum\limits_{i=2}^{r+1}\langle \xi_i,\eta\rangle ^{n-l}\over n-l!}C_l(\xi^2,\dots,\xi^{r+1}, \eta).
\end{gather*}
Because the component of $\phi_{\xi^1}$ along ${\eta^{\otimes
n}\over n!}$ is ${\langle \xi^1,\eta\rangle ^n\over n!}$, this
shows that the component of
\[
:\phi_{\xi^1}\Xi(\phi_{\xi^2},\dots,\phi_{\xi^{r+1}},\lambda\eta):
\]
 along ${\eta^{\otimes n}\over n!}$ is
\begin{gather*}
C_l(\xi^2,\dots,\xi^{r+1}, \eta)\sum_{n_1+n_2=n}{\langle
\sum\limits_{i=2}^{r+1}\xi^i,\eta\rangle^{n_1-l}\over (n_1-l)!} {\langle
\xi^1,\eta\rangle^{n_2}\over n_2!} =
C_l(\xi^2,\dots,\xi^{r+1},\eta){\langle
\sum\limits_{i=1}^{r+1}\xi_i,\eta\rangle^{n-l}\over{(n-l)!}}.
\end{gather*}
The result follows directly.\end{proof}

\begin{definition}The continuous Hochschild cohomology $H^r_{\rm cont}(W_{\infty-},W_{\infty-})$ of the Hida test algebra is the space
 ${\rm Ker}\, \delta^r/{\rm Im}\,\delta^{r-1}$, where the
Hochschild coboundary acts on $L(W_{\infty-}^r,W_{\infty-})$.
\end{definition}

We consider cochains which are {\it finite} sums of
polydif\/ferential operators of degree $(l,m) \in
(\mathbb{N}\times \mathbb{N})$. We call the space of
polydif\/ferential operators $L_{\rm
dif}(W_{\infty-}^r,W_{\infty-})$. By the previous proposition,
$\delta^r$ applies $L_{\rm dif}(W_{\infty-}^r, W_{\infty-})$ into
$L_{\rm dif}(W_{\infty-}^{r+1}, W_{\infty-})$.

\begin{definition}The dif\/ferential Hochschild cohomology $H^r_{\rm dif}(W_{\infty-}, W_{\infty-})$ of the Hida test al\-geb\-ra is the space
${\rm Ker}\, \delta^r/{\rm Im}\,\delta^{r-1}$ where $\delta^r$
acts on $L_{\rm dif}(W_{\infty-}^r, W_{\infty-})$.\end{definition}

We get the main theorem of this work:

\begin{theorem}The differential Hochschild cohomology groups of the Hida test algebra are equal to the continuous
Hochschild cohomology groups of the Hida test algebra.
\end{theorem}

\begin{proof}This comes from the Fock expansion of the previous part and from the following fact: if $\delta\Xi$ is a polydif\/ferential operator
for a continuous cochain $\Xi$, there exists a polydif\/ferential
operator $\Xi_1$ such that $\delta \Xi = \delta \Xi_1$ by
Proposition 1.
\end{proof}

\subsection*{Acknowledgements}

Author thank L.~Accardi and G.~Pinczon for helpful discussions.

\pdfbookmark[1]{References}{ref}
\LastPageEnding

\end{document}